\documentclass[a4paper,12pt]{article}
\usepackage{amssymb,amsthm}

\author{I.~V.~Arzhantsev \quad and \quad D.~A.~Timashev}
\title{Affine embeddings of homogeneous spaces%
\thanks{This work was supported by INTAS--OPEN--97--1570 and
RFBR grant 98--01--00598.}
}
\date{August 29, 2000}

\newcommand{\QQ}{\mathbb{Q}}
\newcommand{\PP}{\mathbb{P}}
\newcommand{\kk}{\Bbbk}

\newcommand{\T}{\mathcal{T}}
\newcommand{\Z}{\mathcal{Z}}
\newcommand{\AF}{\textup{(AF)}}
\newcommand{\g}{\mathfrak{g}}
\newcommand{\h}{\mathfrak{h}}
\newcommand{\sgl}{\mathfrak{sl}}
\newcommand{\z}{\mathfrak{z}}
\newcommand{\e}{\varepsilon}
\newcommand{\rk}{\mathop{\mathrm{rk}}}
\newcommand{\Spec}{\mathop{\mathrm{Spec}}}
\newcommand{\diag}{\mathop{\mathrm{diag}}}
\newcommand{\Hom}{\mathop{\mathrm{Hom}}}
\newcommand{\Aut}{\mathop{\mathrm{Aut}}}
\newcommand{\implies}{\;\Longrightarrow\;}
\newcommand{\codim}{\mathop{\mathrm{codim}}\nolimits}
\newcommand{\trdeg}{\mathop{\mathrm{tr.\,deg}}}
\newcommand{\mod}{\mathop{\mathrm{mod}}\nolimits}

\newtheorem{theorem}{Theorem}
\newtheorem{proposition}{Proposition}
\newtheorem{lemma}{Lemma}
\newtheorem{corollary}{Corollary}
\theoremstyle{definition}
\newtheorem{definition}{Definition}
\newtheorem{example}{Example}

\newtheorem*{Conjecture}{Conjecture}
\theoremstyle{remark}
\newtheorem{remark}{Remark}

\newcounter{property}
\renewcommand{\theproperty}{\textup{(\arabic{property})}}
\newcommand{\property}{\refstepcounter{property}\item}
\newcounter{prooperty}
\renewcommand{\theprooperty}{\textup{(\arabic{prooperty})}}
\newcommand{\prooperty}{\refstepcounter{prooperty}\item}

\makeatletter
\def\keywords#1{{\def\@thefnmark{\relax}\@footnotetext{#1}}}
\let\subjclass\keywords
\makeatother

\begin{document}
\maketitle




\begin{abstract}
Let $G$ be a reductive algebraic group and let $H$ be a reductive
subgroup of $G$. We describe all pairs $(G,H)$ such that
for any affine $G$-variety $X$ with a dense $G$-orbit isomorphic to
$G/H$ the number of $G$-orbits in $X$ is finite. The maximal
number of parameters in families of $G$-orbits in all affine
embeddings of $G/H$ is computed.
\end{abstract}


\section{Introduction.}

Let $G$ be a connected reductive algebraic
group over an algebraically closed
field $\kk$ of characteristic zero and let $H$ be an algebraic subgroup
of $G$. Let us recall that a pointed irreducible
algebraic $G$-variety $X$ is said to be an
\emph{embedding} of the homogeneous space $G/H$ if the base
point of $X$ has the dense orbit and stabilizer~$H$. We
shall denote this by $G/H\hookrightarrow X$.

Let $B$ be a Borel subgroup of $G$. By definition, the
\emph{complexity} $c(X)$ of a $G$-variety $X$ is
the codimension of a generic $B$-orbit in $X$
for the restricted action $B:X$, see~\cite{vin} and~\cite{lv}.
By Rosenlicht's theorem, $c(X)$ is equal to the transcendence
degree of the field $\kk(X)^B$ of rational $B$-invariant
functions on~$X$. A normal $G$-variety $X$ is called
\emph{spherical} if $c(X)=0$, or, equivalently, $\kk(X)^B=\kk$.
A homogeneous space $G/H$ and a subgroup $H\subseteq G$ are said
to be spherical if $G/H$ is a spherical $G$-variety with respect
to the natural $G$-action.

\begin{theorem}[Servedio~\cite{ser}, Luna--Vust~\cite{lv},
Akhiezer~\cite{akh}]
\label{FO}
A homogeneous space $G/H$ is spherical if and only if
each embedding of $G/H$ has finitely many
$G$-orbits.
\end{theorem}

To be more precise, F.~J.~Servedio proved that any affine spherical
variety contains finitely many $G$-orbits, D.~Luna, Th.~Vust
and D.~N.~Akhiezer extended this result to an arbitrary
spherical variety, and D.~N.~Akhiezer constructed
a projective embedding with infinitely many $G$-orbits
for any homogeneous space of positive complexity.

Let us say that
an embedding $G/H\hookrightarrow X$ is \emph{affine} if the variety
$X$ is affine.
In many problems of invariant theory, representation
theory and other branches of mathematics, only affine
embeddings of homogeneous spaces are considered. Hence
for a homogeneous space $G/H$ it is natural to ask:
does there exist an affine embedding $G/H\hookrightarrow X$ with
infinitely many $G$-orbits?

Note that a given homogeneous space $G/H$ admits
an affine embedding if and only if $G/H$ is quasiaffine (as an
algebraic variety), see~\cite[Th.\,1.6]{vp}.
In this situation, the subgroup
$H$ is said to be \emph{observable} in $G$. For a description of
observable subgroups, see \cite{sukh}, \cite[Th.\,4.18]{vp}. By
Matsushima's criterion, $G/H$ is affine iff $H$ is reductive.
(For a simple proof, see~\cite[\S2]{slices}.) In particular,
any reductive subgroup is observable.  In the sequel, we
suppose that $H$ is an observable subgroup of $G$. In this
paper, we are concerned with the following problem:
characterize all quasiaffine homogeneous spaces $G/H$ of a
reductive group $G$ with the property:

\begin{center}
{\AF} \quad
\itshape For any affine embedding $G/H\hookrightarrow X$, the
number of $G$-orbits in $X$ is finite.
\end{center}

\begin{example}\label{sph}
For any spherical quasiaffine homogeneous space,
property {\AF} holds (Theorem~\ref{FO}).
\end{example}

\begin{example}[\cite{po}]\label{SL(2)}
Property {\AF} holds for any homogeneous space
of the group $SL(2)$. In fact, here $\dim X\le 3$, and only a
one-parameter family of one-dimensional orbits can appear in
$X\setminus(G/H)$. But $SL(2)$ contains no two-dimensional
observable subgroups.
\end{example}

\begin{example}\label{max.tor}
Let $T$ be a maximal torus in $G$ and
let $V$ be a finite-dimensional $G$-module. Suppose that a vector
$v\in V$ is $T$-fixed. Then the orbit $Gv$ is closed in $V$,
see \cite{kos}, \cite{lu}. This shows that property {\AF} holds
for any subgroup $H$ such that $T\subseteq H$.
\end{example}

\begin{definition}\label{ac}
An affine homogeneous space $G/H$ is called
\emph{affinely closed} if it admits only one affine embedding
$X=G/H$.
\end{definition}

Homogeneous spaces $G/H$ of Example~\ref{max.tor} are affinely
closed. Denote by $N_G(H)$ the normalizer of $H$ in
$G$. The following theorem generalizes Example~\ref{max.tor}:

\begin{theorem}[Luna~\cite{lu}]\label{fin.ind}
Let $H$ be a reductive subgroup of a reductive group $G$. The
homogeneous space $G/H$ is affinely closed if and only if the
group $N_G(H)/H$ is finite.
\end{theorem}

This theorem provides many examples of
homogeneous spaces with property {\AF}.
Let us note that the complexity of the space $G/T$ can be
arbitrary large, whence property {\AF} cannot be characterized
only in terms of complexity.

In this paper, we show that the union of two conditions---the
sphericity and the finiteness of $N_G(H)/H$---is very close to
characterizing all affine homogeneous spaces of a reductive group
$G$ with property {\AF}. Our main result is

\begin{theorem}\label{main}
For a reductive subgroup $H\subseteq G$, {\AF} holds if and only
if either $N_G(H)/H$ is finite or any extension of $H$ by a
one-dimensional torus in $N_G(H)$ is spherical in~$G$.
\end{theorem}

\begin{corollary}\label{c>1}
For an affine homogeneous space $G/H$ of complexity $>1$, {\AF}
holds iff $G/H$ is affinely closed.
\end{corollary}

\begin{corollary}\label{c=1}
An affine homogeneous space $G/H$ of complexity~$1$ has {\AF}
iff either $N_G(H)/H$ is finite or $\rk N_G(H)/H=1$ and $N_G(H)$
is spherical.
\end{corollary}

The proofs of Theorem~\ref{main} and of its corollaries are
given in Sections~\ref{inf.num}, \ref{proof.main}.

For simple $G$, there is a list of all affine homogeneous spaces
of complexity one~\cite{pan}. We immediately deduce from this
list and Corollary~\ref{c=1} that for simple $G$, there exists
only one series of affine homogeneous spaces of complexity one
that admit affine embeddings with infinitely many $G$-orbits.
Namely, $G=SL(n)$, $n>4$, and $H^0=SL(n-2)\times \kk^*$, where
$SL(n-2)$ is embedded in $SL(n)$ as the stabilizer of the first
two basis vectors $e_1$ and $e_2$ in the tautological
representation of $SL(n)$, and $\kk^*$ acts on $e_1$ and $e_2$
with weights $\alpha_1$ and $\alpha_2$ such that
$\alpha_1+\alpha_2=2-n$, $\alpha_1\ne \alpha_2$, and acts on
$\langle e_3,\dots, e_n\rangle$ by scalar multiplications.

In Section~\ref{verys} we consider very symmetric affine
embeddings ${G/H\hookrightarrow X}$, i.e. affine embeddings whose
group of $G$-equivariant automorphisms $\Aut_G(X)$ contains
the identity component of $\Aut_G(G/H)\cong N_G(H)/H$.
The criterion for finiteness of the number of $G$-orbits in any
such embedding is given (Proposition~\ref{vs}).

The aim of Section~\ref{mod} is to generalize Theorem~\ref{main} following
ideas of~\cite{akh2} and to find the maximal number of parameters
in a continuous family of $G$-orbits over all affine embeddings
of a given affine homogeneous space $G/H$. More precisely,

\begin{definition}\label{modality}
Let $F:X$ be an algebraic group action. The integer
$$
d_F(X)=\min_{x\in X}\codim_XFx=\trdeg\kk(X)^F
$$
is called the {\it generic modality}\/ of the action. The {\it
modality}\/ of $F:X$ is the number $\mod_F X=\max_{Y\subseteq X}
d_F(Y)$, where $Y$ runs through $F$-stable irreducible
subvarieties of $X$.
\end{definition}

Note that $c(X)=d_B(X)$. It was proved by
E.~B.~Vinberg~\cite{vin} that for any $G$-variety $X$ one has
$\mod_B(X)=c(X)$, which means that if we pass from $X$ to a
$B$-stable irreducible subvariety $Y\subseteq X$, then the
number of parameters for $B$-orbits does not increase.

Simple examples show that for $G$ itself
the equality $d_G(X)=\mod_G(X)$ is not true.
This motivates the following

\begin{definition}\label{assa}
With any $G$-variety $X$ we associate the integer
$$
m_G(X)=\max_{X'}\mod_G(X'),
$$
where $X'$ runs through all $G$-varieties birationally $G$-isomorphic
to $X$.
\end{definition}

It is clear that for any subgroup $H\subset G$ the inequality
$m_G(X)\le m_H(X)$ holds. In particular, $m_G(X)\le c(X)$.  The next
theorem shows that $m_G(X)=c(X)$.

\begin{theorem}[Akhiezer~\cite{akh2}]\label{akhi2}
There exists a $G$-variety $X'$ birationally $G$-isomorphic to $X$ such
that $\mod_G(X')=c(X)$.
\end{theorem}

For a homogeneous space $G/H$ we have
$m_G(G/H)=\max_{X}\mod_G(X)$, where $X$ runs through
all embeddings of $G/H$. The affine version of this notion is the following

\begin{definition}\label{bssb}
With any quasiaffine homogeneous space $G/H$ we associate the integer
$$
a_G(G/H)=\max_{X}\mod_G(X),
$$
where $X$ runs through all affine embeddings of $G/H$.
\end{definition}

The following theorem is a direct generalization of Theorem~\ref{main}.

\begin{theorem}\label{gener}
Let $G/H$ be an affine homogeneous space.

(1) If the group $N_G(H)/H$ is finite then $a_G(G/H)=0$;

(2) If $N_G(H)/H$ is infinite then
$$
a_G(G/H)=\max_{H_1}c(G/H_1),
$$
where $H_1$ runs through all non-trivial extensions of $H$ by a one-dimensional
subtorus of $N_G(H)$.
In particular, $a_G(G/H)=$ $c(G/H)$ or $c(G/H)-1$.
\end{theorem}

Applying this theorem to the case $H=\{ e\}$, we obtain

\begin{corollary}\label{vot}
$a_G(G)=\dim U - 1$ if $G$ is semisimple, and $a_G(G)=\dim U$ otherwise,
where $U$ is a maximal unipotent subgroup of $G$.
\end{corollary}


\subsubsection*{Acknowledgements}

We are grateful
to M.~Brion, D.~I.~Panyushev
and E.~B.~Vinberg for useful discussions.

This paper was partially written during our stay
at the Erwin Schr\"odinger International Institute for Mathematical
Physics (Wien, Austria). We wish to
thank this institution for invitation and hospitality.


\subsubsection*{Notation and conventions}

$G$ is a connected reductive group;
\\
$H$ is an observable subgroup of $G$;
\\
$T\subseteq B$ are a maximal torus and a Borel subgroup
of $G$;
\\
$U$ is the maximal unipotent subgroup of $B$;
\\
$N_G(H)$ is the normalizer of $H$ in $G$;
\\
$W(H)$ is the quotient group $N_G(H)/H$;
\\
$\gamma: N_G(H)\to W(H)$ is the quotient homomorphism;
\\
$\kk^*$ is the multiplicative group of non-zero elements of the
base field $\kk$;
\\
$L^0$ is the identity component of an algebraic group $L$;
\\
$Z(L)$ is the center of $L$, $\z(L)$ is its Lie algebra;
\\
$X^L$ is the set of $L$-fixed points in an $L$-variety $X$;
\\
$L_x$ is the isotropy subgroup of  $x\in X$;
\\
$\Xi(G)_+$ is the semigroup of all dominant weights of $G$;
\\
$V_{\mu}$ is an irreducible $G$-module with highest
weight $\mu$;
\\
$\kk[X]$ is the algebra of regular functions and $\kk(X)$ is the
field of rational functions on an algebraic variety~$X$;
\\
$\Spec A$ is the affine variety corresponding to a finitely
generated algebra $A$ without nilpotent elements.

Algebraic groups are denoted by uppercase Latin letters and
their Lie algebras by the respective lowercase Gothic letters.


\section{Embeddings with infinitely many orbits.}
\label{inf.num}

\begin{theorem}\label{basic}
Let $H$ be an observable subgroup in a reductive group $G$.
Suppose that there is a non-trivial one-parameter subgroup
$\lambda: \kk^*\to W(H)$ such that the subgroup
$H_1=\gamma^{-1}(\lambda(\kk^*))$ is not spherical in $G$. Then
there exists an affine embedding $G/H\hookrightarrow X$
with infinitely many $G$-orbits.
\end{theorem}

We shall prove this theorem in the next section. The idea of the
proof is to apply Akhiezer's construction for the non-spherical
homogeneous space $G/H_1$ and to consider the
affine cone over a projective embedding
of $G/H_1$ with infinitely many $G$-orbits.

\begin{proof}[\bf Proof of Corollary~\ref{c>1}]
The assertion follows from Theorem~\ref{fin.ind} and
Theorem~\ref{basic}, which is a part of Theorem~\ref{main} (for
reductive $H$).  Indeed, reductivity of $H$ implies reductivity
of $W(H)$~\cite{lu}. If $W(H)$ is not finite, then it contains
a non-trivial one-parameter subgroup $\lambda(\kk^*)$. For
$H_1=\gamma^{-1}(\lambda(\kk^*))$, we have $c(G/H_1)\ge 1$
whenever $c(G/H)>1$.
\end{proof}

\begin{corollary}\label{center}
Let $G$ be a reductive group with infinite center $Z(G)$ and let
$H$ be an observable subgroup in $G$ that does not contain
$Z(G)^0$. Then  property {\AF} holds for $G/H$ if and only if
$H$ is a spherical subgroup of $G$.
\end{corollary}

\begin{proof}
As $H$ does not contain $Z(G)^0$, there exists
a non-trivial one-parameter subgroup $\lambda(\kk^*)$ in $Z(G)$
with  finite intersection with $H$. The corresponding extension
$H_1$ is spherical iff $H$ is spherical in $G$.
\end{proof}

\begin{corollary}\label{1-ext}
Let $H$ be a connected reductive
subgroup in a reductive group $G$.
Suppose that there exists a reductive non-spherical
subgroup $H_1$ in $G$ such that
$H\subset H_1$ and $\dim H_1=\dim H+1$. Then property {\AF} does
not hold for $G/H$.
\end{corollary}

\begin{proof} Under these assumptions,
there exists a non-trivial one-parameter
subgroup of $H_1$ with finite intersection with $H$ which
normalizes (and even centralizes) $H$.
\end{proof}


\section{Proof of Theorem~\ref{basic}.}

\begin{lemma}\label{fin.ext}
If property {\AF} holds for a homogeneous space $G/H$, then it
holds for any homogeneous space $G/H'$, where $H'$ is an
overgroup of $H$ with $(H')^0=H^0$.
\end{lemma}

\begin{proof} Suppose that there exists an affine embedding
$G/H'\hookrightarrow X$ with infinitely many $G$-orbits.
Consider the morphism $G/H\to G/H'$. It determines an
embedding $\kk[G/H']\subseteq \kk[G/H]$.
Let $A$ be the integral closure of the
subalgebra $\kk[X]\subseteq \kk[G/H']$
in the field of rational functions $\kk(G/H)$.
We have the following commutative diagrams:
$$
\begin{array}{ccccccccc}
A & \hookrightarrow & \kk[G/H] & \hookrightarrow &  \kk(G/H) &
\ \ \ \ \ \ \ & \Spec A & \hookleftarrow & G/H \\
\uparrow &       & \uparrow &  & \uparrow &
\ \ \ \ \ \ \ & \downarrow & & \downarrow \\
\kk[X]       & \hookrightarrow & \kk[G/H'] & \hookrightarrow & \kk(G/H') &
\ \ \ \ \ \ \ & X & \hookleftarrow & G/H'
\end{array}
$$

The affine variety $\Spec A$ with a natural $G$-action
can be considered as an affine embedding
of $G/H$. The embedding $\kk[X]\subseteq A$ defines a finite
(surjective) morphism $\Spec A \to X$ and therefore,
$\Spec A$ contains infinitely many $G$-orbits.
This contradiction completes the proof.
\end{proof}

\begin{remark}\label{fin.res}
The converse statement does not hold.
Indeed, set $G=SL(3)$ and $H=(t, t, t^{-2})\subset T\subset SL(3)$.
We can extend $H$ by a one-parameter subgroup $(t, t^{-1},1)$.
Then $H_1=T$ is not a spherical subgroup in $SL(3)$ and, by
Theorem~\ref{basic}, property {\AF} does not hold here. On the
other hand, one can extend $H$ to $H'$ by a finite noncyclical subgroup
of $W(H)\cong PSL(2)$.
The group $W(H')$ is finite and, by Theorem~\ref{fin.ind},
property {\AF} holds for $G/H'$.
\end{remark}

\begin{lemma}\label{eigenvector}
\begin{enumerate}
\item\label{G_v>H}
Let $H\subseteq G$ be an observable subgroup and $H_1$ be the
extension of $H$ by a one-dimensional torus
$\lambda(\kk^*)\subseteq W(H)$. Then
there exists a finite-di\-men\-sional $G$-module $V$ and an
$H_1$-eigenvector $v\in V$ such that
\begin{list}{\theproperty}{}
\property\label{G<v>} the orbit $G\langle v\rangle$ of the line
$\langle v\rangle$ in the projective space $\PP(V)$ is
isomorphic to $G/H_1$;
\property\label{Hv=v} $H$ fixes $v$;
\property\label{transit} $H_1$ acts transitively on $\kk^*v$.
\end{list}
\item\label{inf.G-orb} If $H_1$ is not spherical in~$G$, then a
couple $(V,v)$ in~\ref{G_v>H} can be chosen so that
\begin{list}{\theproperty}{}
\property\label{inf.orb} the closure of $G\langle v\rangle$ in
$\PP(V)$ contains infinitely many $G$-orbits.
\end{list}
\item\label{G_v=H} If $H$ is reductive, then one can suppose
that $G_v=H$.
\end{enumerate}
\end{lemma}

\begin{proof}
\ref{G_v>H} \
By Chevalley's theorem,
there exists a $G$-module $V'$ and a
vector $v'\in V'$ having property~\ref{G<v>}. Let us
denote by $\chi$ the character of $H$ at~$v'$. Since
$H$ is observable in $G$, every finite-dimensional $H$-module
can be embedded in a finite-dimensional $G$-module~\cite{bbhm}.
In particular, there exists a finite-dimensional $G$-module
$V''$ containing $H$-eigenvectors of character~$-\chi$.
Choose among them a $H_1$-eigenvector $v''$ and put $V=V'\otimes
V''$ and $v=v'\otimes v''$. Properties \ref{G<v>} and~\ref{Hv=v}
are satisfied.

If condition~\ref{transit} also holds, then we are done.
Otherwise, consider any $G$-module $W$ having a vector with
stabilizer~$H$. Take an $H_1$-eigenvector $w\in W^H$ with
nontrivial character, and replace $V$ by $V\otimes W$ and $v$
by $v\otimes w$. Now properties \ref{G<v>}--\ref{transit} are
satisfied.

\ref{inf.G-orb} \
Since $H_1$ is not spherical in $G$, by a result due to
Akhiezer~\cite{akh}, we may choose $(V',v')$ in~\ref{G_v>H} so
that properties \ref{G<v>} and~\ref{inf.orb} are satisfied. Then
we proceed as in~\ref{G_v>H} to obtain the couple $(V ,v)$. The closure
$\overline{G\langle v\rangle}\subseteq\PP(V)$ is contained in the image of
the Segre embedding
$$
  \PP(V')\times\PP(V'')\hookrightarrow\PP(V),
  \qquad\mbox{or}\quad
  \PP(V')\times\PP(V'')\times\PP(W)\hookrightarrow\PP(V),
$$
and projects $G$-equivariantly onto
$\overline{G\langle v'\rangle}\subseteq\PP(V')$. This
implies~\ref{inf.orb} for $(V,v)$.

\ref{G_v=H} \
Let $\omega$ be a fundamental weight of the group $H_1/H$.
Suppose that $H_1/H$ acts at the vector $v$ constructed above by
a character $k\omega$. Since $H_1$ is reductive
(and, in particular, is observable), there exists a $G$-module
$W'$ and an $H_1$-eigenvector $w'\in W'^H$ with weight
$(1-k)\omega$ \cite{bbhm}.
It remains to replace $V$ by $V\otimes W'$ and $v$ by $v\otimes w'$.
\end{proof}

\begin{remark}\label{H=U}
For an arbitrary observable subgroup, statement~\ref{G_v=H} of
Lemma~\ref{eigenvector} does not hold. For example, let $G$ be
the group $SL(3)$ and $H=U$ be a maximal unipotent subgroup
normalized by $T$. Consider the subtorus
$T'=\diag(t^2, t,t^{-3})$
in $T$ as a one-parameter subgroup $\lambda(\kk^*)$.
Any $H$-stable vector in a finite-dimensional $G$-module is a
sum of highest weight vectors. The restriction of any dominant
weight to $T'$ has a non-trivial kernel and the stabilizer of
such a vector contains $H$ as a proper subgroup.
\end{remark}

\begin{proof}[\bf Proof of Theorem~\ref{basic}]
Let $(V,\ v)$ be the couple from Lemma~\ref{eigenvector}. Denote
by $H'$ the stabilizer $G_v$ of the vector $v$. By
\ref{G<v>}--\ref{transit} and since $H_1/H$ is isomorphic to
$\kk^*$, $H'$ is an overgroup of $H$ with $(H')^0=H^0$. By~\ref{transit},
the closure of $Gv$ in $V$ is a cone, so by~\ref{inf.orb} the
property (AF) does not hold for $G/H'$.
Lemma~\ref{fin.ext} completes the proof.
\end{proof}


\section{Proof of Theorem~\ref{main}}
\label{proof.main}

Let $H$ be a reductive subgroup of~$G$. If there exists a
non-spherical extension of $H$ by a one-dimensional torus, then
{\AF} fails for $G/H$ by Theorem~\ref{basic}. To prove the
converse, we begin with the following

\begin{lemma}[{\cite[7.3.1]{val}}]
\label{B-inv}
Let $X$ be an irreducible $G$-variety, and $v$ be a $G$-invariant
valuation of $\kk(X)/\kk$ with residue field~$\kk(v)$. Then
$\kk(v)^B$ is the residue field of the restriction of $v$
to~$\kk(X)^B$.
\end{lemma}

\begin{proof}
For completeness, we give the proof in the case, where $X$ is
affine (the only case we need below). It suffices to prove that
any $B$-invariant element of $\kk(v)$ is the residue class of a
$B$-invariant rational function on~$X$.

For any $f_1,f_2\in\kk(X)$, we shall write $f_1\equiv f_2$
if $v(f_1)=v(f_2)<v(f_1-f_2)$. Such ``congruences'' are
$G$-stable and may be multiplied term by term, as usual
numerical congruences.

Assume $f=p/q$, $p,q\in\kk[X]$, $v(f)=0$, and the residue class
of $f$ belongs to~$\kk(v)^B$. Then $v(p)=v(q)=d$, and
$bf\equiv f$, $\forall b\in B$, i.e.\ $bp\cdot q\equiv p\cdot bq$.

Let $M$ be a complementary $G$-submodule to $\{h\in\kk[X]\mid
v(h)>d\}$ in $\{h\in\kk[X]\mid v(h)\ge d\}$, and $p_0,q_0$ be
the projections of $p,q$ on~$M$. Then $bp_0\cdot q\equiv bp\cdot
q\equiv p\cdot bq\equiv p\cdot bq_0$, $\forall b\in B$. By the
Lie--Kolchin theorem, we may choose finitely many $b_i\in B$,
$\lambda_i\in\kk$ so that $q_1=\sum\lambda_ib_iq_0$ is a
$B$-eigenfunction in $M$ of some weight~$\mu$. Put
$p_1=\sum\lambda_ib_ip_0$. Then $p_1\cdot q\equiv p\cdot q_1$,
whence $p_1/q_1\equiv f\equiv bf\equiv bp_1/\mu(b)q_1$, $\forall
b\in B$. It follows that $bp_1\equiv\mu(b)p_1$, hence
$bp_1=\mu(b)p_1$, because $p_1\in M$. Thus $p_1,q_1$ are
$B$-eigenfunctions of the same weight, and
$f_1=p_1/q_1\in\kk(X)^B$ has the same residue class in $\kk(v)$
as~$f$.
\end{proof}

\begin{definition}[{\cite[\S7]{flats}}]\label{source}
Let $X$ be a normal $G$-variety. A discrete $\QQ$-valued
$G$-invariant valuation of $\kk(X)$ is called \emph{central} if
it vanishes on $\kk(X)^B\setminus\{0\}$. A \emph{source} of $X$
is a non-empty $G$-stable subvariety $Y\subseteq X$ which is the
center of a central valuation of~$\kk(X)$.
\end{definition}

For affine $X$, central valuations are described in a simple
way. Consider the isotypic decomposition
$$
  \kk[X]=\bigoplus_{\mu\in\Xi(X)_{+}}\kk[X]_{\mu},
$$
where the \emph{rank semigroup} $\Xi(X)_{+}\subseteq\Xi(G)_{+}$
is the set of all dominant weights $\mu$ such that
$\kk[X]_{\mu}\ne0$. For any $\lambda,\mu\in\Xi(X)_{+}$, we
have
$$
  \kk[X]_{\lambda}\cdot\kk[X]_{\mu}\subseteq \kk[X]_{\lambda+\mu}
  \oplus \bigoplus_{\alpha\in\T_{\lambda,\mu}(X)}
  \kk[X]_{\lambda+\mu-\alpha},
\leqno(*)
$$
where $\T_{\lambda,\mu}(X)$ is a finite set of positive integral
linear combinations of positive roots, and the inclusion fails
for all proper subsets of~$\T_{\lambda,\mu}(X)$.

Let $\Xi(X)$ be a sublattice spanned by $\Xi(X)_{+}$ in the
weight lattice of~$G$, and $\Xi(X)_{\QQ}=\Xi(X)\otimes\QQ$.
Define the ``cone of tails'' $\T(X)$ to be the convex cone in
$\Xi(X)_{\QQ}$ spanned by the union of all
$\T_{\lambda,\mu}(X)$.

A central valuation $v$ is constant on each $\kk[X]_{\mu}$ and
defines a linear function
$\nu\in\Hom(\Xi(X),\QQ)=\Xi(X)_{\QQ}^{*}$ so that
$\langle\nu,\mu\rangle=v(f)$, $f\in\kk[X]_{\mu}\setminus\{0\}$.
By definition of a valuation, we must have
$\langle\nu,\alpha\rangle\leq0$ for
$\forall\alpha\in\T_{\lambda,\mu}(X)$,
$\lambda,\mu\in\Xi(X)_{+}$. Conversely, each linear function
$\nu\in\Hom(\Xi(X),\QQ)$ which is non-positive on $\T(X)$
defines a central valuation $v$ of $\kk(X)$ by the formula
$$
  v(f)=\min\{\langle\nu,\mu\rangle\mid f_{\mu}\ne0\},
$$
where $f_{\mu}$ is the projection of $f\in\kk[X]$
on~$\kk[X]_{\mu}$.

The valuation $v$ has a center on $X$ iff $\nu$ is non-negative
on $\Xi(X)_{+}$, and the respective source $Y\subseteq X$ is
determined by a $G$-stable ideal
$$
  I(Y)=\bigoplus_{\langle\nu,\mu\rangle>0}\kk[X]_{\mu}
      \mathrel{\;\triangleleft\;}\kk[X]
$$

Central valuations of~$\kk(X)$, identified with respective
linear functions on~$\Xi(X)_{\QQ}$, form a solid convex cone
$\Z(X)\subseteq\Xi(X)_{\QQ}^{*}$, namely the dual cone
to~$-\T(X)$. Knop proved~\cite[9.2]{val}, \cite[7.4]{flats} that
$\Z(X)$ is a fundamental domain for a finite group
$W_X\subset\Aut\Xi(X)$ (the \emph{little Weyl group} of~$X$)
acting on~$\Xi(X)_{\QQ}^{*}$ as a crystallographic reflection
group.

The following lemma is an easy consequence of results of
Knop~\cite{flats}.

\begin{lemma}
\label{1-param}
If $X$ is a normal affine $G$-variety containing a proper
source, then there exists a one-dimensional torus
$S\subseteq\Aut_G(X)$ such that $\kk(X)^B\subseteq\kk(X)^S$.
(Here $\Aut_G(X)$ denotes the group of $G$-equivariant
automorphisms of $X$.)
\end{lemma}

\begin{proof}

If $X$ is as above, then Knop has shown that the algebra
$\kk[X]$ admits a non-trivial $G$-invariant grading, whose
homogeneous components are sums of isotypic components of the
$G$-module $\kk[X]$, see \cite[7.9]{flats} and its proof.
This grading is constructed as follows. Under the above
assumptions, there is a central valuation $v$ of $\kk(X)$ such
that the respective linear function $\nu$ on $\Xi(X)_{\QQ}$ lies
in $\Z(X)\cap-\Z(X)$, hence $\nu$ vanishes on~$\T(X)$. In view
of~$(*)$, this $\nu$ defines a grading of $\kk[X]$ such that
isotypic components $\kk[X]_{\mu}$ are homogeneous of
degree~$\langle\nu,\mu\rangle$.

Let $S\subseteq\Aut_G(X)$ be the one-dimensional torus
corresponding to this grading. Take any $f\in\kk(X)^B$, $f=p/q$,
$p,q\in\kk[X]$. By the Lie--Kolchin theorem, we may choose
finitely many $b_i\in B$, $\lambda_i\in\kk$ so that
$q_0=\sum\lambda_ib_iq$ is a $B$-eigenfunction of some
weight $\mu\in\Xi(X)_{+}$. Then $p_0=\sum\lambda_ib_ip$ is a
$B$-eigenfunction of the same weight, and $f=p_0/q_0$. Since
$p_0,q_0\in\kk[X]_{\mu}$, the torus $S$ acts on them by the same
weight $\langle\nu,\mu\rangle$, thence $f\in\kk(X)^S$. This
shows the inclusion $\kk(X)^B\subseteq\kk(X)^S$.
\end{proof}

\begin{proof}[\bf Proof of Theorem~\ref{main}]
It remains to prove that {\AF} holds for $G/H$ whenever any
extension of $H$ by a one-dimensional torus is spherical. As $H$
is reductive, $W(H)$ is reductive, too. If there exists no
one-parameter extension of $H$ at all, then $W(H)$ is finite and
$G/H$ is affinely closed by Theorem~\ref{fin.ind}. Otherwise
$c(G/H)\le1$. As the spherical case is clear, we may suppose
$c(G/H)=1$.

Let $X$ be an affine embedding of~$G/H$. In order to prove that
$X$ has finitely many $G$-orbits, we may assume that $X$ is
normal. If $X$ contains a proper source, then a one-dimensional
torus $S\subseteq\Aut_G(X)\subseteq\Aut_G(G/H)=W(H)$ provided by
Lemma~\ref{1-param} yields a non-spherical extension of~$H$.
Indeed, if $H_1$ is the preimage of $S$ in $N_G(H)$, then
$\kk(G/H_1)^B=\kk(G/H)^{B\times S}=\kk(X)^{B\times S}=
\kk(X)^B\neq\kk$, since $c(X)=1$. This implies $c(G/H_1)=1$, a
contradiction.

If $X$ contains no proper source, then any proper $G$-stable subvariety
$Y\subset X$ is the center of a non-central $G$-invariant
valuation~$v$. There is an inclusion of residue fields
$\kk(Y)\subseteq\kk(v)\implies \kk(Y)^B\subseteq\kk(v)^B$. By
Lemma~\ref{B-inv}, $\kk(v)^B$ is the residue field of the
restriction of $v$ to $\kk(G/H)^B$, which is the field of
rational functions in one variable. As $v$ is non-central,
$\kk(Y)^B=\kk(v)^B=\kk$, thence $Y$ is spherical. It follows
that $X$ has finitely many orbits.  (Otherwise, a one-parameter
family of $G$-orbits provides a non-spherical $G$-subvariety.)
\end{proof}

\begin{proof}[\bf Proof of Corollary~\ref{c=1}]
The reductive group $W(H)$ acts on $\kk(G/H)^B$, which is the
field of rational functions on a projective line. If the kernel
of this action has positive dimension, then it contains a
one-dimensional torus extending $H$ to a non-spherical subgroup.

Otherwise, either $W(H)$ is finite or $\rk W(H)=1$ and each
subtorus of $W(H)$ has a dense orbit on the projective line. The
corollary follows.
\end{proof}

\begin{remark}
In the proof of Theorem~\ref{main}, we have used reductivity
of $H$ only in the following assertion:
\begin{quote}
If $W(H)$ contains no subtori, then it is finite, and $G/H$ is
affinely closed.
\end{quote}
In fact, we need this assertion only if $c(G/H)>1$.
Theorem~\ref{main} holds for quasiaffine $G/H$ of
complexity~$\le1$.

One might hope that the situation described in the above
assertion never occurs for non-reductive~$H$, i.e.\ that $W(H)$
always contains a subtorus. Unfortunately, $W(H)^0$ may be a
non-trivial unipotent group, as the following example shows.
\end{remark}

\begin{example}\label{uni}
Let $e$ be a regular nilpotent in the Lie algebra $\sgl(3)$,
$G=SL(3)\times SL(3)$, and $H$ be the two-dimensional unipotent
subgroup with the Lie algebra generated by $(e, e^2)$ and $(e^2,
e)$.  Then the Lie algebra of the normalizer of $H$ is the
linear span of $(e, 0), \ (e^2, 0), \ (0, e)$ and $(0, e^2)$.
Hence the group $W(H)^0$ is two-dimensional and unipotent.
\end{example}
(Another example was suggested by E.A.Tevelev.)

\medskip

We are not able to characterize quasiaffine, but not affine,
homogeneous spaces with the property (AF).

\medskip

In this context we would like to formulate the following

\begin{Conjecture}
If $H\subseteq G$ is observable, but not reductive, then $W(H)$
is infinite.
\end{Conjecture}


\section{Very symmetric embeddings.}
\label{verys}

The group of $G$-equivariant automorphisms of a homogeneous
space $G/H$ is isomorphic to $W(H)$. (The action $W(H):G/H$
is induced by the action $N_G(H):G/H$ by right multiplication.)
Let $G/H\hookrightarrow X$ be an affine embedding. The group
$\Aut_G X$ of $G$-equivariant automorphisms of $X$ is a subgroup
of $W(H)$.

\begin{definition}\label{vs-def}
An embedding $G/H\hookrightarrow X$ is said to be
\emph{very symmetric} if $W(H)^0\subseteq\Aut_G X$.
\end{definition}

Any spherical affine variety is very symmetric.
In fact, for a spherical homogeneous space $G/H$, any isotypic
component $\kk[G/H]_{\mu}$ of the $G$-algebra $\kk[G/H]$
is an irreducible $G$-module (see~\cite{ser} or~\cite[Th.2]{vk}), and $W(H)$
acts on $\kk[G/H]_{\mu}$ by scalar multiplications.
This shows that any $G$-invariant subalgebra in $\kk[G/H]$
is $W(H)$-invariant, too.

In the case of affine $SL(2)/\{e\}$-embeddings,
only the embedding $X=SL(2)$ is very symmetric; in all other
cases the group $\Aut_{SL(2)} X$ is isomorphic to a Borel
subgroup in $SL(2)$, see~\cite[III.4.8, Satz 1]{kr}.
More generally, if $X$ is an affine embedding of the homogeneous
space $G/\{e\}$, then $X$ is very symmetric if and only if
the action $G:X$ can be extended to an action of the group
$G\times G$ with an open orbit isomorphic to $(G\times G)/H$,
where $H$ is the diagonal in $G\times G$.
Hence $X$ can be considered as an affine $(G\times G)/H$-embedding.
Theorem~\ref{fin.ind} implies that if $G$ is a semisimple
group, then $X=(G\times G)/H$, for other proofs see~\cite{wat}
and~\cite[Prop.\,1]{semi}.

If $G$ is a reductive
group, then the set of all very symmetric embeddings
of the homogeneous space $G/\{e\}$ is
exactly the set of all affine algebraic monoids with $G$
as the group of units~\cite{semi}.
Thus very symmetric embeddings have a natural
characterization in the variety of all affine
$G/\{e\}$-embeddings.
The classification of reductive algebraic monoids is obtained
in~\cite{semi} and~\cite{rit}.

Put $\widehat{G}=G\times W(H)^0$,
$N=\gamma^{-1}(W(H)^0)$, and
$\widehat{H}=\{(n,nH)\mid n\in N\}$ (the ``diagonal''
embedding of~$N$). Any very symmetric
affine embedding of $G/H$ may be considered as an embedding of
$\widehat{G}/\widehat{H}$.

\begin{proposition}\label{vs.inf}
Under assumptions of Theorem~\ref{basic}, if $\lambda(\kk^{*})$
is central in $W(H)^0$, then there exists a very symmetric affine
embedding $G/H\hookrightarrow X$ with infinitely many $G$-orbits.
\end{proposition}

\begin{proof}
We follow the proof of Theorem~\ref{basic}. Put
$\widehat{H}_1=\widehat{H}\cdot\lambda(\kk^{*})$; then
$\widehat{H}_1\cap G=H_1$.
We modify the proof of Lemma~\ref{eigenvector}\ref{inf.G-orb} to
obtain a $\widehat{G}$-module $V$ and an
$\widehat{H}_1$-eigenvector $v\in V$ such that
$\widehat{G}_{\langle v\rangle}=\widehat{H}_1$, $\widehat{G}_v$
is a finite extension of~$\widehat{H}$, and
$\overline{\widehat{G}\langle v\rangle}\subseteq\PP(V)$ contains
infinitely many $G$- (not $\widehat{G}$-) orbits. Arguing as in
the proof of Theorem~\ref{basic}, we see that the closure $X$ of
$Gv=\widehat{G}v\subseteq V$ is $\widehat{G}$-stable and has
infinitely many $G$-orbits, q.e.d. (Observe that $\widehat{G}$
may be not reductive, but Lemma~\ref{fin.ext}, required in the
proof, does not use the reductivity assumption.)

To construct such a couple $(V,v)$, it suffices, in the notation
of Lemma~\ref{eigenvector}, to construct a $\widehat{G}$-module
$V'$ and a vector $v'\in V'$ such that $G\langle
v'\rangle=\widehat{G}\langle v'\rangle\cong
\widehat{G}/\widehat{H}_1$ and $\overline{G\langle v'\rangle}$
has infinitely many $G$-orbits. Then we proceed as in
Lemma~\ref{eigenvector}\ref{G_v>H}, replacing $G$
by~$\widehat{G}$. (Note that the reductivity of $G$ is not
essential in Lemma~\ref{eigenvector}\ref{G_v>H}.) It remains to
construct a couple $(V',v')$. For this purpose, we refine
Akhiezer's construction~\cite{akh}.

By assumption, $c(G/H_1)>0$, hence there exists a character
$\xi:H_1\to\kk^*$ such that for the associated line bundle
$L_{\xi}$ on $G/H_1$, the multiplicity of a certain simple
$G$-module $V_{\mu}$ in $H^0(G/H_1,L_{\xi})$ is greater than
one~~\cite[Th.\,1]{vk}.

The group $W(H)^0$ acts on
$H^0(G/H_1,L_{\xi})$ and on the isotypic component
$E=H^0(G/H_1,L_{\xi})_{\mu}$ by $G$-module automorphisms.

Take a $\widehat{G}$-module $M$ and a vector $m\in M$ such that
$\widehat{G}_{\langle m\rangle}=\widehat{H}_1$. Let $Y$ be the
closure of $\widehat{G}\langle m\rangle=G\langle m\rangle$
in~$\PP(M)$. The natural rational map $f:Y\dasharrow\PP(E^{*})$
is $\widehat{G}$-equivariant.

Consider a decomposition
$E=E_0\oplus\dots\oplus E_k$ into
irreducible $G$-submodules and fix isomorphisms
$\psi_i:V_{\mu}\to E_i$. Choose a basis
$\{\e_0,\dots,\e_m\}$ of $T$-eigenvectors with weights
$\mu_0=\mu,\mu_1,\dots,\mu_m$ in~$V_{\mu}$, and put
$\e^{(i)}_j=\psi_i(\e_j)$. In projective coordinates,
$$
  f(gH_1)=[\e^{(0)}_0(gH_1):\dots:\e^{(k)}_0(gH_1):\dots:
           \e^{(0)}_m(gH_1):\dots:\e^{(k)}_m(gH_1)]
$$

The closure $Z$ of the graph of $f$ in
$Y\times\PP(E^{*})$ is $\widehat{G}$-stable. We claim that $Z$ contains
infinitely many $G$-orbits.
To prove it, take a strictly dominant one-parameter subgroup
$\delta:\kk^{*}\to T$. If all $\e^{(i)}_0(gH_1)\neq0$, then
\begin{eqnarray*}
f(\delta(t)gH_1)&=&[\,\cdots\,:\,t^{-\langle\mu_j,\delta\rangle}
\e^{(i)}_j(gH_1)\,:\,\cdots\,]\\
&=&[\,\cdots\,:\,t^{\langle\mu_0-\mu_j,\delta\rangle}
\e^{(i)}_j(gH_1)\,:\,\cdots\,]\\
&\longrightarrow&
[\e^{(0)}_0(gH_1):\dots:\e^{(k)}_0(gH_1):\dots:0:\dots:0]
\end{eqnarray*}
as $t\to0$, because $\mu_0-\mu_j$ is a positive linear
combination of positive roots for all $j>0$.

We may identify
$E^{*}$ with $V_{\mu}^{*}\otimes\kk^{k+1}$ and consider the
Segre embedding
$\PP(V_{\mu}^{*})\times\PP^k\hookrightarrow\PP(E^{*})$. Then
$\lim_{t\to0}\delta(t)f(gH_1)=(\langle\e_0^{*}\rangle,p)\in
\PP(V_{\mu}^{*})\times\PP^k$, where $\{\e_j^{*}\}$ is the dual
basis to~$\{\e_j\}$, and
$p=[\e^{(0)}_0(gH_1):\dots:\e^{(k)}_0(gH_1)]\in\PP^k$.

As the sections $\e^{(0)}_0,\dots,\e^{(k)}_0$ are linearly
independent on~$G/H_1$, $\overline{f(Y)}$ intersects infinitely
many closed disjoint $G$-stable subvarieties
$\PP(V_{\mu}^{*})\times\{p\}\hookrightarrow\PP(E^{*})$,
$p\in\PP^k$. This proves the claim, because $Z$ projects
$G$-equivariantly onto~$\overline{f(Y)}$.

Finally, a $\widehat{G}$-module $V'=M\otimes E^{*}$ and a vector
$v'=m\otimes e$ such that $f(\langle m\rangle)=\langle e\rangle$
are the desired, because $\overline{G\langle v'\rangle}\cong Z$.

The proof is complete.
\end{proof}

Now we are interested in the following problem: when does
any very symmetric affine embedding of a homogeneous space $G/H$
have finitely many $G$-orbits? The example of
$SL(3)/\{e\}$-embeddings shows that the latter property is not
equivalent to {\AF}.

\begin{proposition}\label{vs}
Let $H$ be a reductive subgroup in a reductive
group $G$. Every very symmetric affine embedding of $G/H$ has
finitely many $G$-orbits iff either {\AF} holds or $W(H)^0$ is
semisimple. In the second case, there is only one very symmetric
affine embedding, namely $X=G/H$.
\end{proposition}

\begin{proof}
The Lie algebra of $N_{\widehat{G}}(\widehat{H})$ equals
$\widehat{\h}+\widehat{\z}$, where $\widehat{\z}$ is the
centralizer of $\widehat{H}$ in~$\widehat{\g}$. We have
$\widehat{\z}=\z(N)\oplus\z(W(H)^0)$, and
$\z(N)\cong\z(H)\oplus\z(W(H)^0)$.

If $W(H)^0$ is semisimple, then $\widehat{\z}\subseteq\h
\subseteq\widehat{\h}$, and $N_{\widehat{G}}(\widehat{H})$ is
finite. Theorem~\ref{fin.ind} implies the assertion for this
case.

Now suppose that $W(H)^0$ is not semisimple. If there exists a
non-spherical extension of $H$ by a one-dimensional torus
$S\subseteq Z(N)$, then by Proposition~\ref{vs.inf},
there exists a very symmetric affine embedding of $G/H$ with
infinitely many $G$-orbits.

Finally, suppose that any extension of $H$ by a one-dimensional
torus in $Z(N)$ is spherical. Then $c(G/H)\le1$. As
the spherical case is clear, we may assume that $c(G/H)=1$.

The connected kernel $W_0$ of the action $W(H):\kk(G/H)^B$ acts
on isotypic components of $\kk[G/H]$ by scalar multiplications.
Whence $W_0$ is diagonalizable and central in $W(H)$. By
assumption, $W_0=\{e\}$. Hence $W(H)^0$ is a one-dimensional
torus acting on $\kk(G/H)^B$ with finite kernel. By
Corollary~\ref{c=1}, {\AF} holds for~$G/H$. The proof is
complete.
\end{proof}

\section{Affine embeddings and modality}
\label{mod}

We begin this section with the generalization of Lemma~\ref{eigenvector}.

\begin{lemma}\label{general}
Let $H\subseteq G$ be an observable subgroup and $H_1$ be the
extension of $H$ by a one-dimensional torus
$\lambda(\kk^*)\subseteq W(H)$. Then
there exists a finite-di\-men\-sional $G$-module $V$ and an
$H_1$-eigenvector $v\in V$ such that
\begin{list}{\theprooperty}{}
\prooperty\label{G<v>1} the orbit $G\langle v\rangle$ of the line
$\langle v\rangle$ in the projective space $\PP(V)$ is
isomorphic to $G/H_1$;
\prooperty\label{Hv=v1} $H$ fixes $v$;
\prooperty\label{transit1} $H_1$ acts transitively on $\kk^*v$;
\prooperty\label{inf.G-orb1} $\mod_G(\overline{G\langle v\rangle})=c(G/H_1)$.
\end{list}
\end{lemma}

\begin{proof}
We use exactly the same arguments as in the proof of
Lemma~\ref{eigenvector} replacing an embedding of
$G/H_1$ with infinitely many orbits
from~\cite{akh} by an embedding of $G$-modality $c(G/H_1)$
constructed in~\cite{akh2}.
\end{proof}

\begin{lemma}\label{ggff}
In the notation of Lemma~\ref{general},
$$
c(G/H)\geq
a_G(G/H)\geq c(G/H_1)\geq c(G/H)-1
$$
In particular, $a_G(G/H)=c(G/H)$ or $c(G/H)-1$.
\end{lemma}

\begin{proof}
Clearly, $a_G(G/H)\le c(G/H)$. Taking an affine cone over the projective
embedding constructed in Lemma~\ref{general}, one obtains an
affine embedding of $G/H'$ of modality $\geq c(G/H_1)$, where
$H'=G_v$ is a finite extension of $H$. Using the
construction from the proof of Lemma~\ref{fin.ext}, we get an affine
embedding of $G/H$ of modality $\geq c(G/H_1)$. The obvious
inequality $c(G/H_1)\ge c(G/H)-1$ completes the proof.
\end{proof}

\begin{proof}[Proof of Theorem~\ref{gener}]
Statement (1) follows from
Theorem~\ref{fin.ind}. To prove (2), we can use Lemma~\ref{ggff}.
If there exists a one-dimensional torus in $N_G(H)$ such that the
extension $H\subset H_1$ is non-trivial and $c(G/H)=c(G/H_1)$,
then there exists an affine embedding of $G/H$ of modality $c(G/H)$.

Conversely, suppose that $G/H\hookrightarrow X$ is an affine embedding
of modality $c(G/H)$. We need to find a one-dimensional subtorus
$S\subseteq W(H)$ such that for the extended subgroup $H_1$ we
will have $c(G/H_1)=c(G/H)$. By the definition of modality, there
exists a proper $G$-invariant subvariety $Y\subset X$, such that
the codimension of a generic $G$-orbit in $Y$ is $c(G/H)$.
Therefore, $c(Y)=c(G/H)$.

Consider a $G$-invariant valuation
$v$ of $\kk(X)$ with the center $Y$. For the residue field
$\kk(v)$ we have $\trdeg\kk(v)^B\geq\trdeg\kk(Y)^B$, hence
$\trdeg\kk(v)^B=\trdeg\kk(X)^B$. If the restriction of
$v$ to $\kk(X)^B$ is not trivial, then by Lemma~\ref{B-inv},
$\trdeg\kk(v)^B<\trdeg\kk(X)^B$, a contradiction. Thus $v$ is
central, and $Y$ is a source of $X$. A one-dimensional
subtorus $S\subseteq\Aut_G(X)\subseteq\Aut_G(G/H)=W(H)$ provided
by Lemma~\ref{1-param} yields the extension of $H$ of the same
complexity.
\end{proof}

\begin{proof}[Proof of Corollary~\ref{vot}]
If $G$ is not semisimple, then for a central one-dimensional
subtorus $T_1$ one has $c(G/T_1)=c(G)=\dim U$. If $G$ is
semisimple, then for any one-dimensional subtorus $T_1\subset G$
there exists a Borel subgroup $B$ which does not contain~$T_1$,
and there is a $B$-orbit on $G/T_1$ of dimension $\dim B$. This
implies $c(G/T_1)=c(G)-1$.
\end{proof}



\begin{thebibliography}{BBHM}

\bibitem[Akh1]{akh}
D.~N.~Akhiezer, Actions with a finite number of orbits,
Funkts. Analiz i ego Prilog. {\bf 19} (1985), no.~1, 1--5
(in Russian); English transl.: Func. Anal. Appl. {\bf 19}
(1985), no.~1, 1--4.

\bibitem[Akh2]{akh2}
D.~N.~Akhiezer, On modality and complexity of reductive group actions,
Uspekhi Mat. Nauk {\bf 43:2} (1988), 129--130
(in Russian); English transl.: Russian Math. Surveys {\bf 43:2}
(1988), 157--158.


\bibitem[BBHM]{bbhm}
A.~Bialynicki-Birula, G.~Hochschild and G.~D.~Mostow, Extensions of
representations of algebraic linear groups, Amer. J. Math. {\bf 85}
(1963), 131--144.

\bibitem[KV]{vk}
B.~N.~Kimel'feld and E.~B.~Vinberg, Homogeneous domains on flag
manifolds and spherical subgroups of semisimple Lie groups,
Funkts. Analiz i ego Prilog. {\bf 12} (1978), no.~3, 12--19
(in Russian); English transl.: Func. Anal. Appl. {\bf 12}
(1978), no.~3, 168--174.

\bibitem[Kn1]{val}
F.~Knop, \"Uber Bewertungen, welche unter einer reduktiven
Gruppe invariant sind, Math. Ann. {\bf 295} (1993), 333--363.

\bibitem[Kn2]{flats}
F.~Knop, The asymptotic behavior of invariant collective motion,
Inv. Math. {\bf 116} (1994), 309--328.

\bibitem[Kos]{kos}
B.~Kostant, Lie group representations on polynomial rings,
Amer. J. Math. {\bf 85} (1963), 327--404.

\bibitem[Kr]{kr}
H.~Kraft, Geometrische Methoden in der Invariantentheorie,
Vieweg Verlag, Braunschweig 1985.

\bibitem[Lu1]{slices}
D.~Luna, Slices \'etales,
Mem. Soc. Math. France {\bf 33} (1973), 81--105.

\bibitem[Lu2]{lu}
D.~Luna, Adh\'erences d'orbite et invariants,
Inv. Math. {\bf 29} (1975), 231--238.

\bibitem[LV]{lv}
D.~Luna and Th.~Vust, Plongements d'espaces homog\`enes,
Comment. Math. Helv. {\bf 58} (1983), 186--245.

\bibitem[Pan]{pan}
D.~I.~Panyushev, Complexity of quasiaffine homogeneous $G$-varieties,
$t$-decompositions, and affine homogeneous spaces
of complexity 1, Advances in Soviet Math., vol.~8,
ed. by E.~B.~Vinberg (1992), 151--166.

\bibitem[Po]{po}
V.~L.~Popov, Quasihomogeneous affine algebraic varieties of the
group $SL(2)$, Izv. Akad. Nauk SSSR, Ser. Mat. {\bf 37} (1973),
no.~4, 792--832 (in Russian); English transl.: Math.
USSR-Izv. {\bf 7} (1973), no.~4, 793--831.

\bibitem[PV]{vp}
V.~L.~Popov and E.~B.~Vinberg, Invariant theory, Itogi Nauki i Tekhniki,
Sovr. Probl. Mat. Fund. Napravl.,
vol.~55, VINITI, Moscow 1989, pp.~137--309 (in Russian); English
transl.: Algebraic Geometry IV, Encyclopaedia of Math. Sciences,
vol.~55, Springer-Verlag, Berlin 1994, pp.~123--278.

\bibitem[Rit]{rit}
A.~Rittatore, Algebraic monoids and group embeddings,
Transformation Groups {\bf 3} (1998), no.~4, 375--396.

\bibitem[Ser]{ser}
F.~J.~Servedio, Prehomogeneous vector spaces and varieties,
Trans. Amer. Math. Soc. {\bf 176} (1973), 421--444.

\bibitem[Su]{sukh}
A.~A.~Sukhanov, Description of the observable subgroups of
linear algebraic groups,
Mat. Sbornik {\bf 137} (1988), no.~1, 90--102 (in Russian);
English transl.: Math. USSR-Sb. {\bf 65} (1990), no.~1, 97--108.

\bibitem[Vi1]{vin}
E.~B.~Vinberg, Complexity of actions of reductive groups,
Funkts. Analiz i ego Prilog. {\bf 20} (1986), no.~1, 1--13
(in Russian); English transl.: Func. Anal. Appl. {\bf 20}
(1986), no.~1, 1--11.

\bibitem[Vi2]{semi}
E.~B.~Vinberg, On reductive algebraic semigroups,
In: Lie Groups and Lie Algebras, E.~B.~Dynkin Seminar,
Amer. Math. Soc. Transl. (2)  {\bf 169} (1995), 145--182.

\bibitem[Wat]{wat}
W.~C.~Waterhouse, The unit group of affine algebraic monoids,
Proc. Amer. Math. Soc. {\bf 85} (1982), 506--508.

\end{thebibliography}
\end{document}